\documentclass[12pt]{article}

\usepackage[margin=1in]{geometry}  
\usepackage{graphicx}              
\usepackage{amsmath}               
\usepackage{amsfonts}              
\usepackage{amsthm}                
\usepackage{tikz}
\usepackage{hyperref}
\usepackage{float}

\theoremstyle{definition}

\title{Efficient construction by ruler and compass of a Pentadecagon}
\author{Arjeh Kurzweil, Erez Sheiner}

\begin{document}

\maketitle

\begin{abstract}
In this short paper we show that with a small change of the common ruler and compass construction of the regular pentadecagon, we can produce more regular polygons .\\
\end{abstract}

\section{Introduction}
Ruler and compass constructions date back to ancient Greek mathematicians, and was studied by Euclid in his famous Elements books from 300 BC (\cite{S},\cite{T}).\\
Using the ruler one can draw a straight line between two points. Using the compass and two given points one can draw a circle around one of the points, with radius which is the distance between theses two points.\\

Among other problems, Euclid studied the problem of constructing regular (equilateral and equiangular \cite{S}) polygons. Only in 1837, Wantzel has proved that a regular polygon can be constructed if and only if the number of its sides is the product of a power of 2 and distinct Fermat primes (\cite{W}).\\

\begin{figure}[H]\centering
\includegraphics[scale=4]{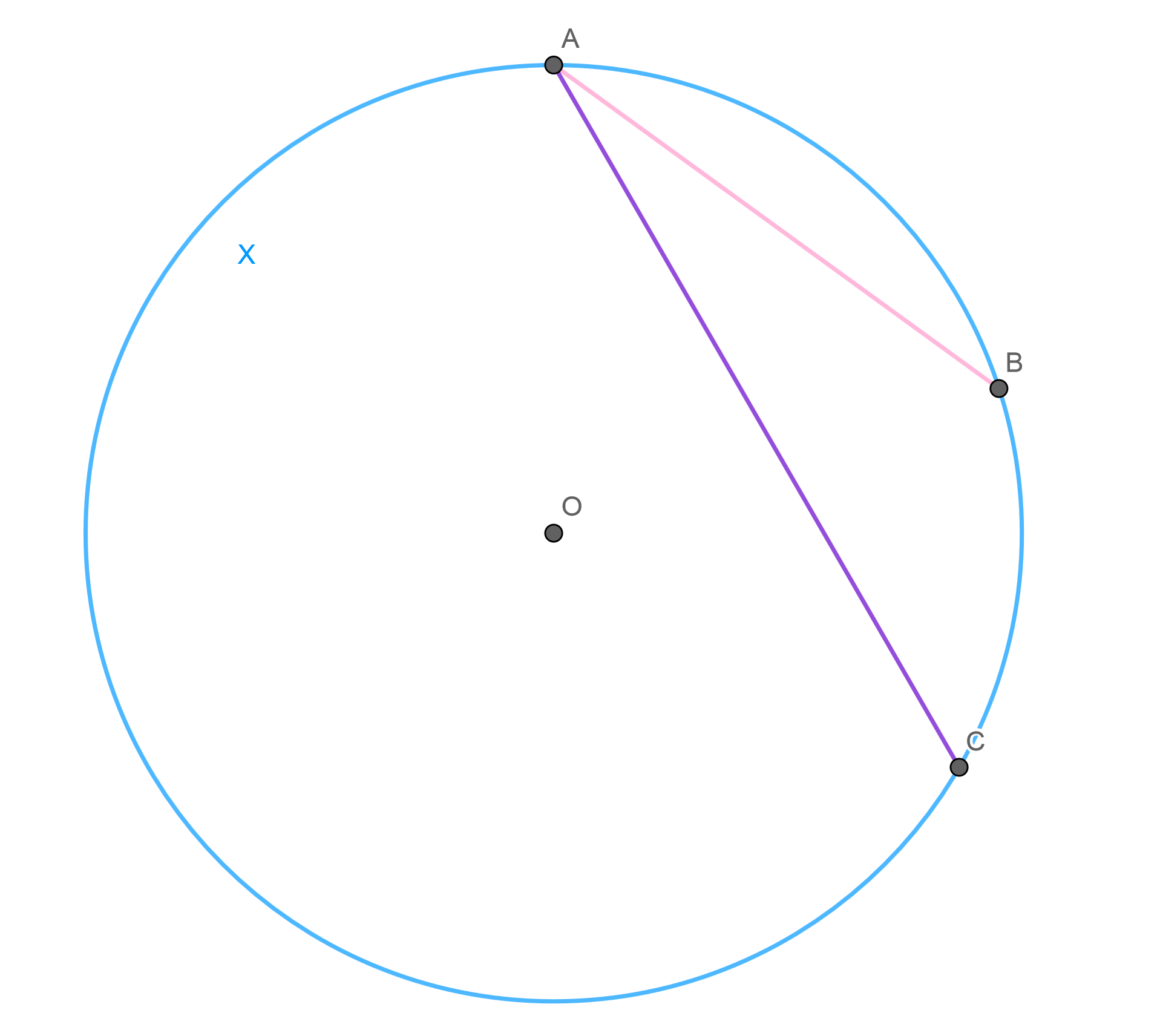}
\caption{}
\label{fig0}
\end{figure}

In Euclids Elements fourth book, problem XVI, he describes how to construct a regular pentadecagon (quindecagon \cite{S}) inscribed in a given circle. Euclid used the construction of a triangle and a pentagon from a given point $A$ (see figure \ref{fig0}). The arc $AB$ is a fifth of the whole circumference, and the arc $AC$ is a third. Therefore the arc $BC$ is two fifteenths of the whole, and by bisecting the angle $\angle BOC$ we obtain the angle of the pentadecagon.\\

Euclid main interest here was the ability to construct the polygon. In this short paper we show that while constructing the regular pentadecagon we can obtain many regular polygons on the way.\\

\section{The construction}

We start with the common construction. Given two arbitrary points $O$ and $A$, we use a compass to construct the circle $X$ with radius $OA$ around the center $O$. We will construct regular polygons inscribed in this circle.\\
Mark the intersection of the line $OA$ with the circle as $B$.\\

We draw circles with radius $AB$ around $A$ and $B$, mark their intersections by $C$ and $D$ and draw the segment between them.\\

\begin{figure}[H]\centering
\includegraphics[scale=2]{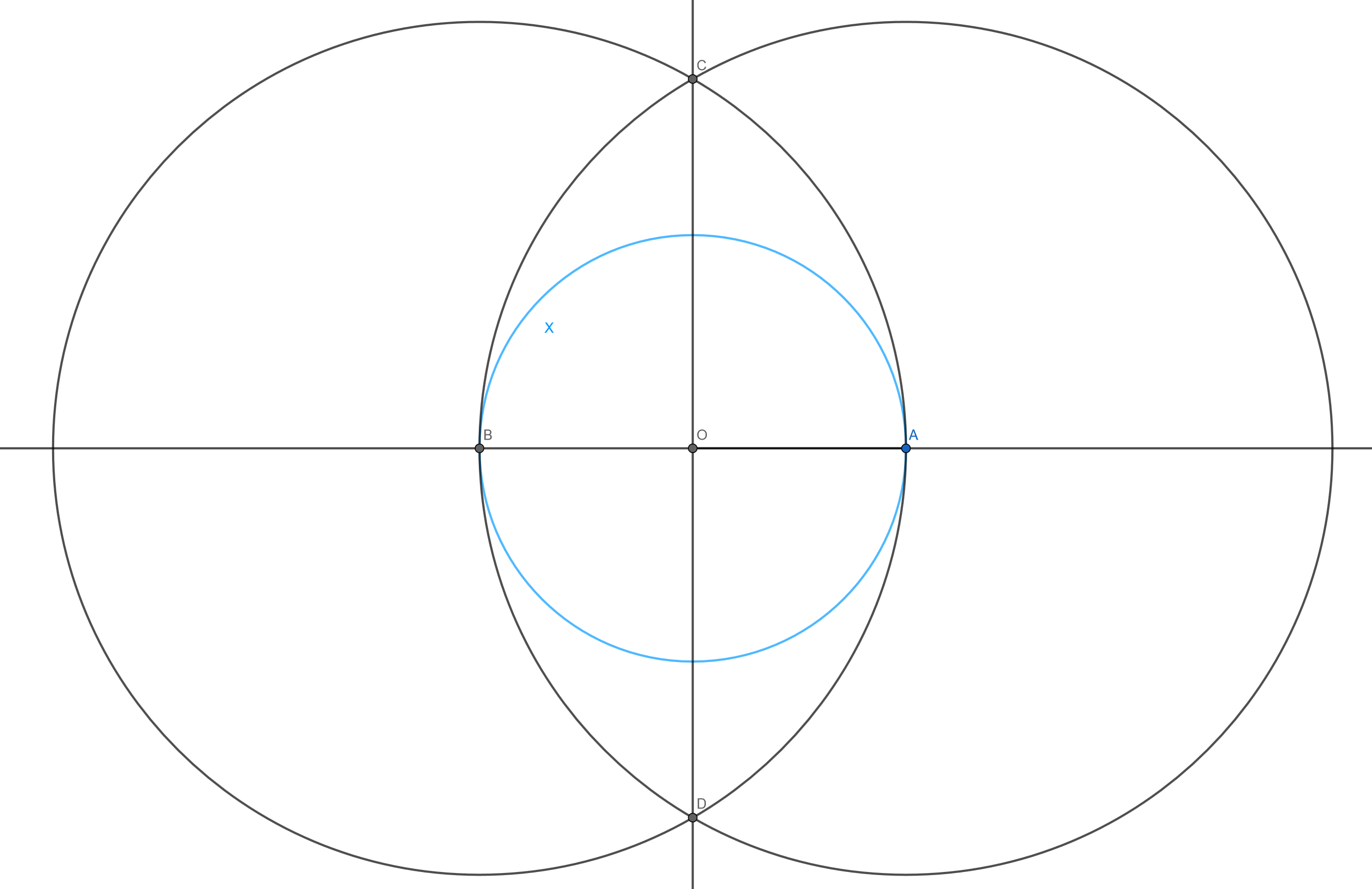}
\caption{}
\label{fig1}
\end{figure}

\begin{figure}[H]\centering
\includegraphics[scale=3]{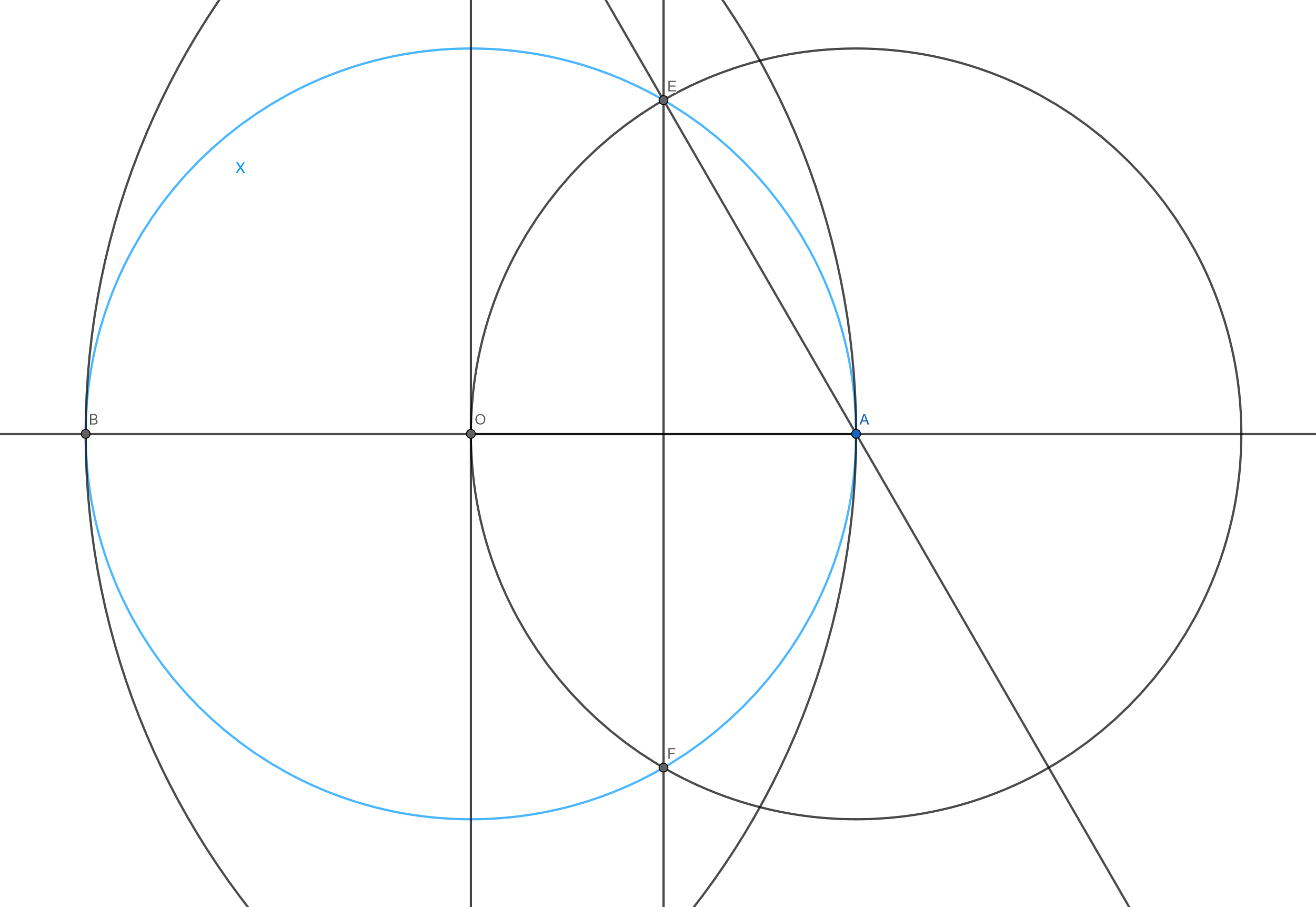}
\caption{}
\label{fig2}
\end{figure}

Next we draw the circle with radius $OA$ around $A$, and mark its intersections with the main circle $X$ by $E$ and $D$ (Figure \ref{fig2}). Now we have that the lengths $OA,OE,AE$ are all equal (the distance between the center of a circle and any point on the circle is equal to the radius), and so we obtain the angles $\angle EOA = 60^\circ$ and $\angle EOF=120^\circ$.\\

Therefore $EF$ is an edge of a regular triangle inscribed in $X$, and $EA$ is an edge of a regular hexagon.\\

\begin{figure}[H]\centering
\includegraphics[scale=2.5]{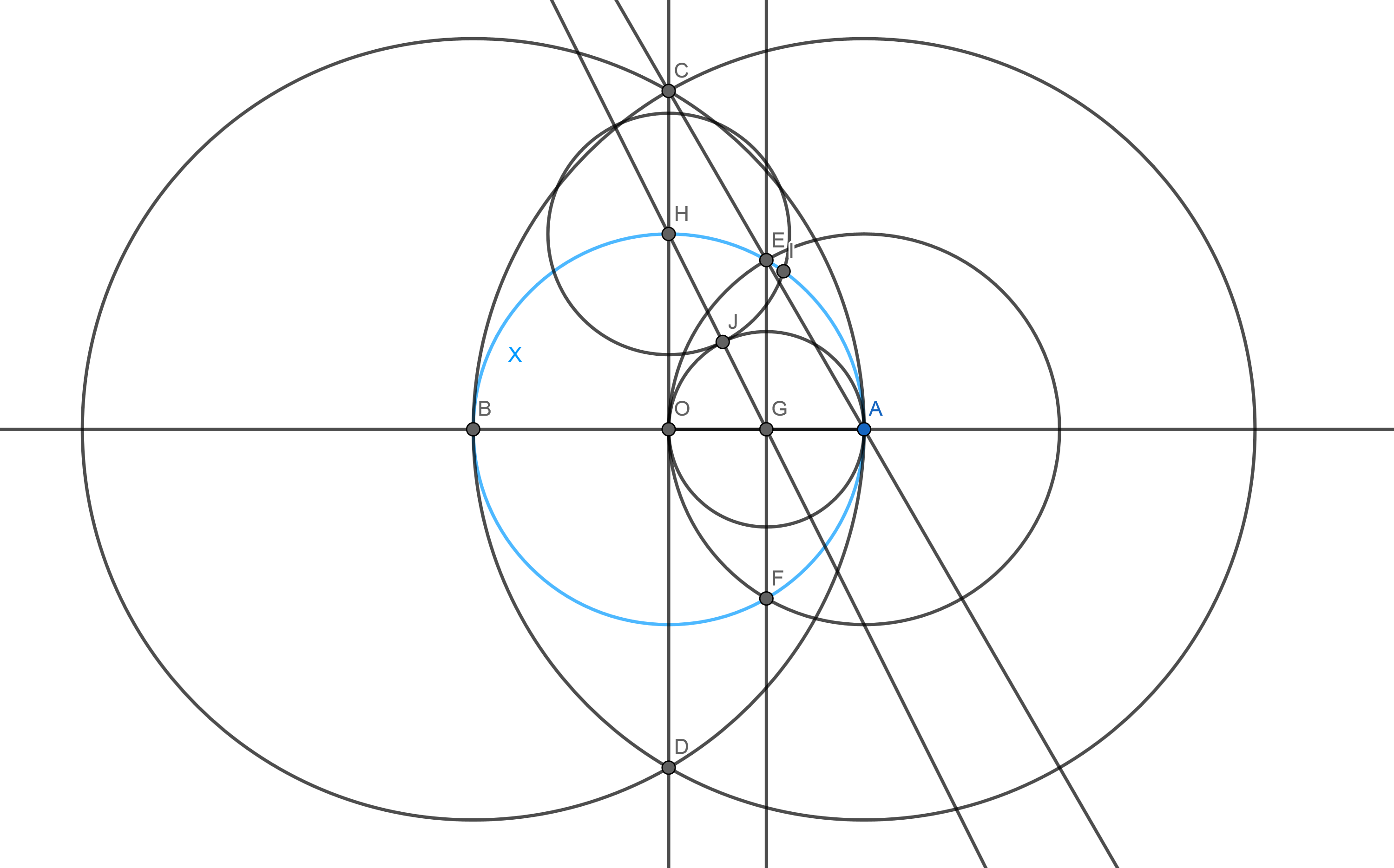}
\caption{}
\label{fig3}
\end{figure}

Mark the intersection of $EF$ and $OA$ by $G$, and the intersection of $CD$ and $X$ by $H$ (Figure \ref{fig3}). Clearly $\angle HOA=90^\circ$ and therefore $HA$ is an edge of an inscribed square.\\

Since $\angle EOA=60^\circ$ we get $\angle HOE=30^\circ$ so $HE$ is an edge of an inscribed regular dodecagon.\\

Draw a circle of radius $OG$ around $G$, and mark its intersection with the line $HG$ by $J$. Draw a circle of radius $HJ$ around $H$ and mark its intersection with $X$ by $I$.\\

\begin{figure}[H]\centering
\includegraphics[scale=2]{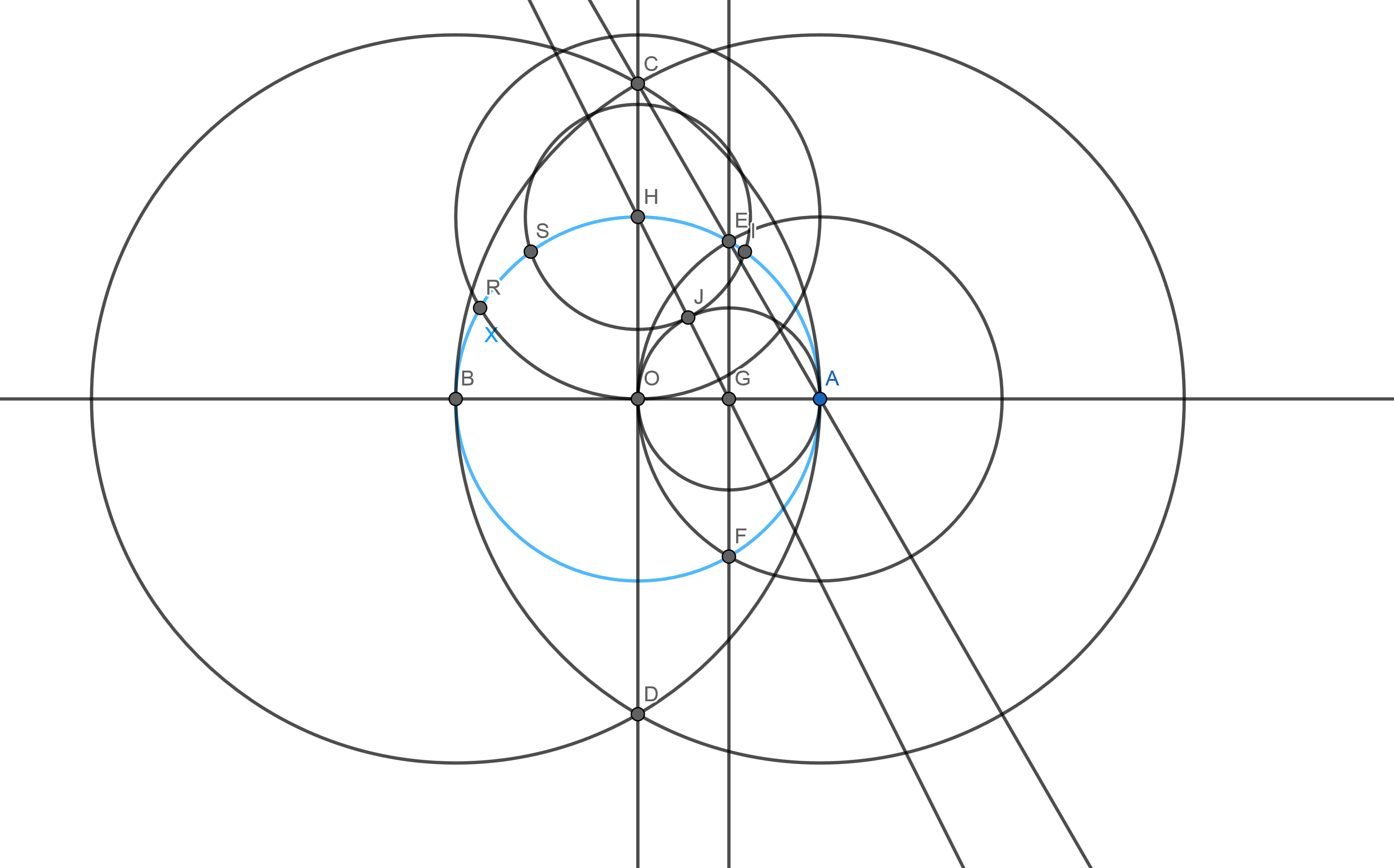}
\caption{}
\label{fig5}
\end{figure}

At this point the common construction of the pentadecagon will continue by creating a circle of radius $OH$ around $H$, intersecting the circle $X$ at $R$. This way we obtain $\angle ROS=24^\circ$, thus $RS$ is the edge of the regular pentadecagon (see figure \ref{fig5}). However, we do not need this extra construction.\\

Since $OG=JG=\frac{1}{2}$ and $HG=\sqrt{1+\frac{1}{4}}$, then $HJ=\frac{\sqrt{5}-1}{2}$ which is exactly the length of an edge of a regular decagon. Therefore $HI$ (which is the same length as $HJ$) is an edge of an inscribed decagon, and $\angle HOI=36^\circ$. We can easily construct the edge of a pentagon by doubling angle $\angle HOI$ and creating the point $K$. \\

Together with the fact that $\angle EOA=60^\circ$ and $\angle HOA=90^\circ$ we conclude that \linebreak $\angle HOE=30^\circ$ and $\angle EOI=6^\circ$. We now know that $EI$ is an edge of a regular $60-gon$, and we easily construct the edges of a regular $30-gon$ and pentadecagon by creating points $L,M,N$ such that $EI=IL=LM=MN$ (see figure \ref{fig4}).\\

We will sum up the edges we constructed in Table \ref{tab1}, illustrated in Figure \ref{fig4}.\\
\begin{table}[H]
  \begin{center}
    \begin{tabular}{l|c|r} 
      $\textbf{n-gon}$ &\textbf{Edge} & \textbf{Angle}\\
      \hline
      3 & $EF$ & $120^\circ$\\
      4 & $HA$ & $90^\circ$\\
      5 & $HK$ & $72^\circ$\\
      6 & $EA$ & $60^\circ$\\
      10 & $HI$ & $36^\circ$\\
      12 & $HE$ & $30^\circ$\\
      15 & $EN$ & $24^\circ$\\
      20 & $IN$ & $18^\circ$\\
      30 & $LN$ & $12^\circ$\\
      60 & $EI$ & $6^\circ$\\

    \end{tabular}
  \end{center}
  \caption{}
    \label{tab1}
\end{table}

\begin{figure}[H]\centering\label{fig4}
\includegraphics[scale=6.3]{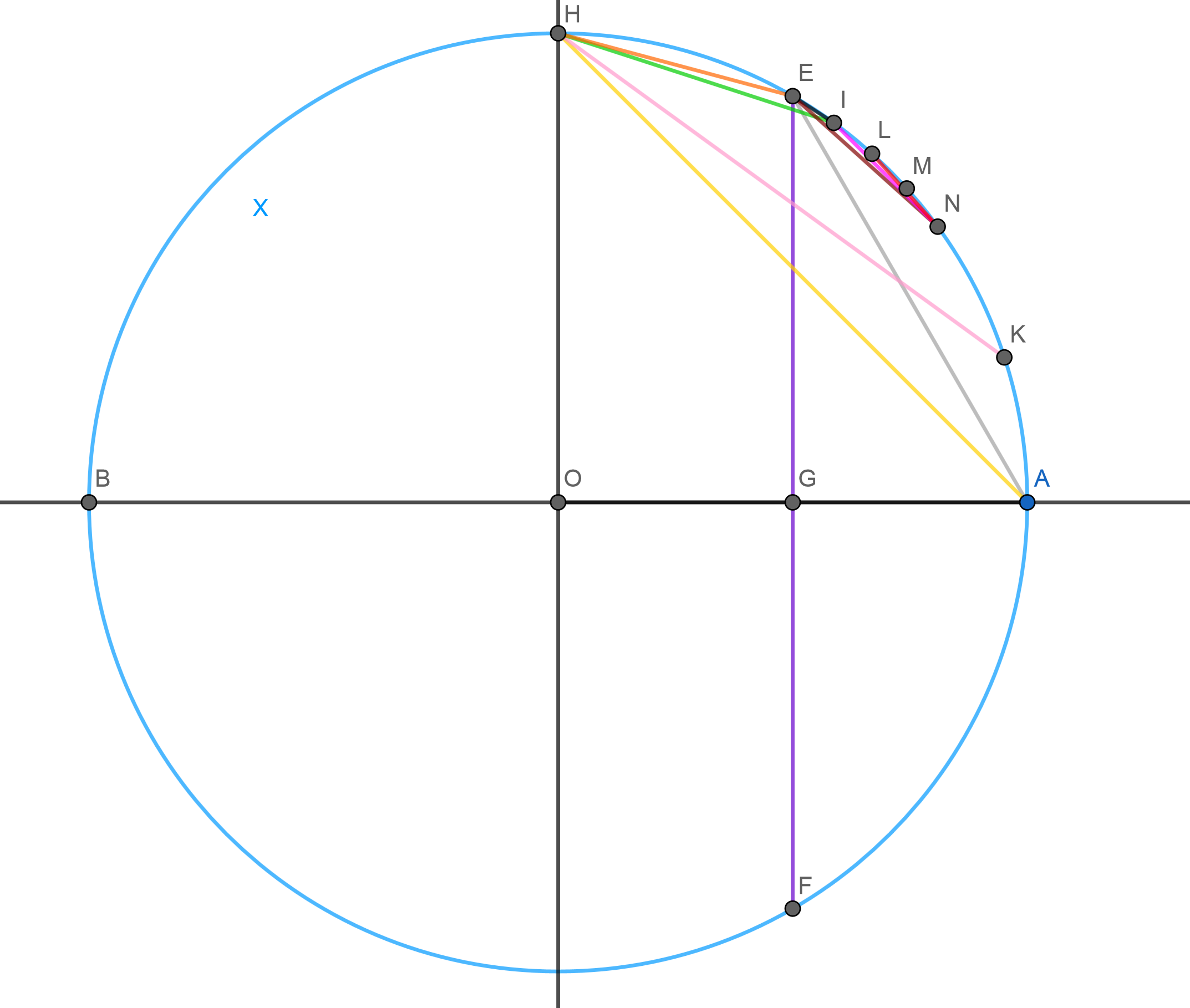}
\caption{}
\label{fig4}
\end{figure}


\clearpage

\begin{thebibliography}{10}	


\bibitem{S} R. Simson, \emph{The elements of Euclid}, Desilver, Thomas and Company (1838).	
\bibitem{T} G. Toussaint, \emph{A new look at Euclid’s second proposition}, The Mathematical Intelligencer 15.3 (1993), 12-24.
\bibitem{W} P.L. Wantzel, \emph{Recherches sur les moyens de reconnaître si un problème de géométrie peut se résoudre avec la règle et le compas}, Journal de Mathématiques pures et appliquées 2.1 (1837), 366-372.





\end{thebibliography}
\end{document}